\documentclass[11pt,a4paper,oneside]{article}
\usepackage{amsmath}
\usepackage{amsfonts}
\usepackage{amsthm}
\usepackage{color}

\newtheorem{thm}{Theorem}

\begin{document}
\title{Homological Properties of the Algebra of Compact Operators on a Banach Space}
\author{G.A. WILLIS}

\maketitle
\paragraph{Abstract.} The conditions on a Banach space, $E$, under which the algebra, $\mathcal{K}(E)$, of compact operators on $E$ is right flat or homologically unital are investigated. These homological properties are related to factorization in the algebra and it is shown that, for $\mathcal{K}(E)$, they are closely associated with the approximation property for $E$. The class of spaces, $E$, such that $\mathcal{K}(E)$ is known to be right flat and homologically unital is extended to include spaces which do not have the bounded compact approximation property.

\section*{Introduction}
\paragraph*{}This note is concerned with the flatness and H-unitality of certain algebras of operators on Banach spaces. The algebra of all operators on the Banach space $E$ will be denoted $\mathcal{B}(E)$. Most of the discussion will deal with the norm closure of the ideal of finite rank operators on $E$, which will be denoted $\mathcal{F}(E)$. The ideal of compact operators on $E$, which might not equal $\mathcal{F}(E)$ when $E$ does not have the approximation property, will be denoted $\mathcal{K}(E)$. 

\paragraph*{}Flatness and H-unitality have been studied by M. Wodzicki, see [Wol, 2 and 3], for a number of algebras arising in functional analysis. They are important in the proof of Karoubi's conjecture of the equality of the algebraic and topological K-theory groups of stable $C^*$-algebras, see \cite{SW1,SW2}. Theorems 1 and 2 below give a partial solution to a question, posed in [Wo2], concerning the H-unitality of $\mathcal{K}(E)$.


The present note is essentially the preprint \cite{Wi3} with the bibliography updated and some other minor changes. The original version was  withdrawn from submission with the aim of  revising it under the guidance of M.~Wodzicki. That intention was not carried through however because the author was involved in other projects and M.~Wodzicki came to see it as part of a much larger and deeper project on categories of Banach spaces \cite{Wo4}. Following further discussion with M.~Wodzicki, it has been decided to publicise the answer to the original question in this version pending what is hoped to be a more complete development in a subsequent publication. 

\section*{Some Homological Algebra}
Let $A$ be an algebra and $X$ be a right $A$-module. Then $X$ is said to be a $flat$ module if, for every exact sequence of left $A$-modules,  
$$0 \longrightarrow M' \stackrel{v}\longrightarrow M \stackrel{w}\longrightarrow M'' \longrightarrow 0, $$
the sequence
$$0\longrightarrow X\otimes_A M' \stackrel{I\otimes v}\longrightarrow X\otimes_A M \stackrel{I\otimes w}\longrightarrow X\otimes_A M'' \longrightarrow 0 $$
is exact, see [Bo, I.2.3. definition 2]. The characterisation of flat modules which it will be most convenient to use here is that given in [Bo, I.2.11, corollary 1]. The right $A$-module $X$ is flat if and only if , whenever $\{ e_i\}_{i\in \mathbb{I}}$ and $\{ b_i\}_{i\in \mathbb{I}}$ are elements of $X$ and $A$ respectively such that $\sum_{i \in \mathbb{I}} e_ib_i =0$, there are elements of $X$,$\{x_j\}_{j\in \mathbb{J}}$ and elements of $A$, $\{a_{j,i} \}_{(j \in \mathbb{J}, i \in \mathbb{I})}$ such that $\sum_{i \in \mathbb{I}} a_{j,i} b_i =0$ for each $j$ and $e_i = \sum_{j \in \mathbb{J}} x_j a_{j,i}$ for each $i$.

\paragraph*{} The algebra $A$ is said to be $right$ $universally$ $flat$ if, whenever $A$ is embedded as a right ideal in an algebra $B$, $A$ is a flat right $B$-module. It is easy to check that, if $A$ has a left unit, then it is right universally flat. There is a corresponding notion of flatness for left $A$-modules and left universally flat algebras.

\paragraph*{}For an algebra, $A$, let $A^{\otimes n}$ denote the $n$-fold tensor product $A\otimes \cdot\cdot\cdot \otimes A$ and define a map $\delta_n :A^{\otimes n} \to A^{\otimes (n-1)}$ by
$$\delta_n(a_1\otimes \cdot\cdot\cdot \otimes a_n) = \sum_{j=1}^{n-1} (-1)^{j-1} a_1 \otimes \cdot\cdot\cdot \otimes a_j a_{j+1} \otimes \cdot\cdot\cdot \otimes a_n.$$
Then $A$ is said to be $homologically$ $unital$, abbreviated $H$-$unital$ if the sequence
\begin{equation} 0 \longleftarrow A \stackrel{\delta_2}\longleftarrow A\otimes A \stackrel{\delta_3}\longleftarrow A\otimes A \otimes A \stackrel{\delta_4}\longleftarrow \cdot\cdot\cdot \end{equation}
is exact.

\paragraph*{}The notion of H-unitality was introduced by M. Wodzicki in [Wo1, Wo2]. It is shown in those papers that H-unitality of $A$ is equivalent to $A$ having the excision property in various homology theories, see \cite[Proposition 2]{Wo1} and [Wo2, theorem 3.1].

\paragraph*{}If $A$ is right (or left) universally flat, then it is H-unital. To see this, suppose that $A$ is right universally flat and let $A^\#$ denote the algebra obtained by adjoining a unit, $1$, to $A$. Then $A$ is a right ideal in $A^\#$. Given $v=\sum_{j\in\mathbb{J}} a_{j,1}\otimes \cdot\cdot\cdot \otimes a_{j,n}$ in $A^{\otimes n}$, choose, for each $k=3,..., n$ a basis, $\{b_{1,k},..., b_{m_k, k}\}$, for span$\{a_{j,k-1}, a_{j,k-1}a_{j,k}, a_{j,k}:j\in \mathbb{J}\}$. Let $\mathbb{B} = \{b_{j_3,3}\otimes \cdot\cdot\cdot \otimes b_{j_n,n}: j_k = 1,2, ..., m_k; k= 3,4, ..., n \}$. Then $\mathbb{B}$ is a linearly independent set in $A^{\otimes(n-2)}$ and 
$$\delta_n(v) = \sum_{j\in \mathbb{J}} \sum_{k=1}^{n-1} (-1)^{k-1} a_{j,1} \otimes \cdot\cdot\cdot \otimes a_{j,k} a_{j, k+1} \otimes \cdot\cdot\cdot \otimes a_{j,n} = \sum_{b\in \mathbb{B}} (\sum_{j\in\mathbb{J}} a_{j,1}c_{j,b} )\otimes b,$$
where each $c_{j,b}$ belongs to $A^\#$ for each $j\in \mathbb{J}$ and each $b$ in $\mathbb{B}$. Hence, if we now suppose that $\delta_n (v) = 0$, then $\sum_{j\in \mathbb{J}} a_{j,1}c_{j,b} = 0$ for each $b$ in $\mathbb{B}$. This is equivalent to $\sum_{j\in \mathbb{J}} a_{j,1}\otimes \tilde{c_j} = 0$ in $A\otimes_{A^\#}F$, where $F$ denotes the direct sum of $|\mathbb{B}|$ copies of $A^\#$ and $\tilde{c_j} = (c_{j,b})_{b \in \mathbb{B}}$. Since $A$ is right universally flat, it is a flat right $A^\#$-module and so, by [Bo, I.2.3, proposition 1], $A$ is $F$-flat. Hence, by [Bo, I.2.11, proposition 13], there are a family $\{x_p\}_{p\in\mathbb{P}}$ of elements of $A$ and a family $\{f_{p,j}\}_{(p\in \mathbb{P}, j\in\mathbb{J})}$ of elements of $A^\#$ such that
\begin{equation}\sum_{j\in\mathbb{J}} f_{p,j} \tilde{c_j} =0 \quad \text{for all}\,\, p \in \mathbb{P}, \, \text{and}  \end{equation} 
\begin{equation} a_{j,1} = \sum_{p\in\mathbb{P}} x_p f_{p,j} \quad \text{for all}\,\, j\in \mathbb{J}. \end{equation}

Now it follows from (2) that 
$$\delta _n(\sum_{j\in \mathbb{J}}f_{p,j}\otimes a_{j,2}\otimes \cdot \cdot \cdot \otimes a_{j,n}) = \sum_{b\in \mathbb{B}} (\sum_{j\in \mathbb{J}} f_{p,j} c_{j,b})\otimes b = 0,$$
for each $p$ in $\mathbb{P}$, hence, if we define $w$ in $A\otimes A^\# \otimes A^{\otimes(n-1)}$ by
$$w = \sum_{p \in \mathbb{P}}\sum_{j \in \mathbb{J}} x_p \otimes f_{p,j} \otimes a_{j,2} \otimes \cdot \cdot \cdot \otimes a_{j,n}, $$
then, by (3),
\begin{align*}\delta_{n+1}(w)  & = \sum_{j \in \mathbb{J}}(\sum_{p \in \mathbb{P}}x_p f_{p,j}) \otimes a_{j,2} \otimes \cdot \cdot \cdot \otimes a_{j,n}\\ & - \sum_{p \in \mathbb{P}} x_p \otimes \delta_n (\sum_{j \in \mathbb{J}} f_{p,j} \otimes a_{j,2} \otimes \cdot\cdot\cdot \otimes a_{j,n}) \\ & = \sum_{j \in \mathbb{J}} a_{j,1} \otimes a_{j,2} \cdot\cdot \cdot \otimes a_{j,n} = v. \end{align*}

Now write $f_{p,j} = c1 + f'_{p,j}$, where $c$ is a scalar and each $f'_{p,j}$ belongs to $A$ and put 
$$ w' = \sum_{p \in \mathbb{P}}\sum_{j \in \mathbb{J}} x_p \otimes f'_{p,j} \otimes a_{j,2} \otimes \cdot \cdot \cdot \otimes a_{j,n}.$$

Then $w'$ belongs to $A^{\otimes (n-1)}$ and, as is easily checked, $\delta_{n+1}(w') = v$. Therefore the sequence (1) is exact at $A^{\otimes n}$ and so $A$ is H-unital.

\paragraph*{} In section 8 of [Wo8] it is asked for which Banach spaces $E$ is it the case that $\mathcal{K}(E)$ is H-unital. This questions is partially answered in [Wo2] and [Wo3]. The stronger result is in [Wo3] where a certain factorization property, ($\Phi$), is given which, if satisfied by an algebra, implies that the algebra is right universally flat, and hence H-unital. It follows from the proof of Cohen's factorization theorem, [BD, Theorem 11.10], that a Banach algebra with a bounded left approximate identity satisfies ($\Phi$). Now $\mathcal{K}(E)$ has a bounded left approximate identity if and only if $E$ has the bounded compact approximation property, see [Di, Theorem 2.6], and so $\mathcal{K}(E)$ is right universally flat if $E$ has this property, [Wo3, Theorem 8(f)].  
\paragraph*{} It is shown below, in theorem 1, that $\mathcal{K}(E)$ can be right universally flat when $E$ does not have the bounded compact approximation property. However, in the circumstances considered $E$ must have the approximation property, in which case $\mathcal{F}(E) = \mathcal{K}(E)$. Consequently all of the following discussion will be concerned with when $\mathcal{F}(E)$ is right universally flat. Note that the bounded compact approximation property does not imply the approximation property, see [Wi1], and so there are spaces to which the result of Wodzicki applies which are not covered by Theorem 1.
\paragraph*{} For the new cases covered by theorem 1, $E$ does not have the bounded approximation property and so $\mathcal{F}(E)$ does not have a bounded left approximate identity [Di, Theorem 2.6]. It is not clear whether $\mathcal{F}(E)$ satisfies the factorization property ($\Phi$) in these cases. Right universal flatness is proved by another method which at one point requires an automatic continuity lemma.
\section*{An Automatic Continuity Lemma}
\paragraph*{}The lemma concerns right multipliers, where a right multiplier on an algebra, $A$, is a linear map $T:A\to A$ such that $T(ab) = aT(b)$ for all $a$ and $b$ in $A$. Note that, if an algebra $A$ is a right ideal in a larger algebra $B$, then for each $b$ in $B$, the map $a\mapsto ab:A\to A$ is a right multiplier on $A$ and it is in this context that we will be applying the lemma.

\newtheorem*{lemma}{Lemma}
\begin{lemma}
Let $E$ be a Banach space and $T$ be a right multiplier on $\mathcal{F}(E)$. Then $T$ is continuous.
\end{lemma}

\it Proof.  \rm Following [Sinc], section 1, let
\begin{align*}\mathfrak{G}(T) = \{a\in \mathcal{F}(E) & : \text{there are}\,\, a_n \in \mathcal{F}(E), n=1,2,...\,\text{s.t.} \lim_{n\to\infty} ||a_n|| =0 \\ & \,\, \text{and} \lim_{n\to \infty} ||a-T(a_n)|| =0 \} \end{align*}
be the separating space of $T$ and $\mathfrak{J}_T = \{ a\in \mathcal{F}(E) \, :\, a\mathfrak{G}(T) = (0) \}$ be the continuity ideal of $T$. Then, since $T$ is a right multiplier, $\mathfrak{G}(T)$ is a closed left ideal and $\mathfrak{J}_T$ is a closed two-sided ideal in $\mathcal{F}(E)$.
\paragraph*{}The lemma is obvious if $E$ is finite dimensional and so we may suppose that $E$ is infinite dimensional and choose sequences $\{e_n\}_{n=1}^\infty$ in $E$ and $\{e_n^*\}_{n=1}^\infty$ in $E^*$ such that $\langle e_i^*, e_j \rangle = \delta_{i,j}$. Then $\{e_n \otimes e_n^* \}_{n=1}^\infty$ is a sequence of mutually orthogonal rank one projections on $E$. It follows from an application of lemma 1.6 in [Sinc] that there is an $n$ such that $e_n \otimes e_n^*$ belongs to $\mathfrak{J}_T$. Since $\mathcal{F}(E)$ is topologically simple it follows that $\mathfrak{J}_T = \mathcal{F}(E)$, whence the separating space of $T$ is zero. Therefore, by lemma 1.2 in [Sinc], $T$ is continuous.
\paragraph*{}For another proof of the lemma in the cases in which we will apply it, see corollary 4 below. Although it will not be required in what follows, note that the set of all multipliers on $A$ is an algebra under composition and we have the following

\newtheorem{corollary}{Corollary}
\begin{corollary}
The algebra of right multipliers of $\mathcal{F}(E)$ is isomorphic to $\mathcal{B}(E^*)$.
\end{corollary}

\it Proof. \rm For each operator $T$ in $\mathcal{B}(E^*)$ define a right multiplier, $\phi(T)$, on $\mathcal{F}(E)$ by $\phi(T)(e\otimes e^*) = e\otimes Te^*$, for rank one operators $e\otimes e^*$ and then extending to all of $\mathcal{F}(E)$ by linearity and continuity. Then $\phi$ is clearly an injective algebra homomorphism from $\mathfrak{B}(E^*)$ into the algebra of right multipliers on $\mathcal{F}(E)$. 
\paragraph*{}To show that $\phi$ is a surjection, let $p=e\otimes e^*$ be a rank one projection on $E$ and $M$ be a right multiplier on $\mathcal{F}(E)$. Then, since $M(pa) = pM(a)$ for all $a$ in $\mathcal{F}(E)$, it follows that $M$ leaves the space of $e\otimes E^*$ invariant. This subspace is isomorphic to $E^*$ and so the restriction of $M$ to the subspace determines an operator, $T$, on $E^*$ which is continuous by the lemma. It is easily checked that $\phi(T) = M$.
\paragraph*{}Similar arguments to the above show that all left multipliers on $\mathcal{F}(E)$ are continuous and that the algebra of all left multipliers on $\mathcal{F}(E)$ is isomorphic to $\mathcal{B}(E)$.

\section*{Right Universal Flatness of $\mathcal{F}(E)$}
As has already has been pointed out, the factorization property $\Phi$ implies right universal flatness. On the other hand, if an algebra $A$ is H-unital, then exactness of the sequence (1) at $A$ implies that $A$ is the span of its products. It is thus clear that factorization is closely connected with the homological properties we are interested in. For this reason, the spaces $C_p$, $1\leq p \leq \infty$, introduced by W. B. Johnson in [Jo] will be important\footnote{ These spaces have also played an important role in the work of A.~Blanco, in which he shows that the algebras of approximable operators considered in Theorem~\ref{thm:one} are weakly amenable whether $E$ has the approximation property or not, \cite{Bl1,Bl2}.}. 
These spaces have the property that, for any Banach space $E$, every operator in $\mathcal{F}(E)$ factors through $C_p$, see theorem 1 in [Jo].
\paragraph*{}In order to describe these spaces, first define a distance, $d(X,Y)$, between isomorphic Banach spaces $X$ and $Y$ by $d(X,Y)=\inf \{||T||\,||T^{-1}||\, : T\, \text{is an isomorphism from}\, X \text{to}\, Y \}$. Now let $\{ G_i \}_{i=1}^\infty$ be a sequence of finite dimensional Banach spaces such that:
\\
i) for each finite dimensional space $E$ each $\epsilon > 0$ there is an $i$ such that $d(E,G_i)< 1+\epsilon$; 
\\
and
\\
ii) for each $i$ there are infinitely many $j\neq i$ such that $G_j$ is isometric to $G_i$.
\\
Next define, for $1\leq p \leq \infty$, $C_p = (\bigoplus_{i=1}^\infty G_i)_p$, where the subscript $p$ indicates the $\ell_p$ direct sum. Up to isomorphism, the space $C_p$ is independent of the choice of sequence of finite dimensional spaces $\{ G_i \}_{i=1}^\infty$, see [Sg, definition 12.3].
\paragraph*{}Many theorems about the spaces $C_p$ describe what are called in [Sing] their \lq universal complement properties\rq , for examples see [Jo, theorem 4] and [Sing, theorem 12.4 and proposition 13.13]. The content of all these results is that, if a Banach space $E$ has some variant of approximation property, then $E\oplus C_p$ has a better approximation property. The theorem we now prove is of a similar kind.

\begin{thm}
\label{thm:one}
For each Banach space $E$ the following are equivalent:
\\
\\
$(1)$ $E$ has the approximation property;
\\
$(2)$ $\mathcal{F}(E\oplus C_1)$ is right flat; \,\, and
\\
$(3)$ $\mathcal{F}(E\oplus C_1)$ is right universally flat.
\end{thm}
    
\it Proof. \rm We begin by showing that (1) implies (3). Suppose that $E$ has the approximation property. Then $E\oplus C_1$ has the approximation property.
\paragraph*{}Let $\mathcal{F}(E\oplus C_1)$ be a right ideal in the algebra $B$ and let $T_1,T_2, ...\, , T_n$ in $\mathcal{F}(E\oplus C_1)$ and $b_1,b_2, ...\, , b_n$ in $B$ be such that $\sum_{k=1}^n T_kb_k =0$. Define the subspace, $Y$ of $\mathcal{F}(E\oplus C_1)^{\oplus n}$ by
$$Y = \{(S_1.S_2, ...\, , S_n): S_k \in \mathcal{F}(E\oplus C_1) \,\, \text{and}\,\, \sum_{k=1}^n S_kb_k =0  \}.$$
Then $(T_1,T_2, ...\, , T_n)$ belongs to $Y$ and since, by the lemma, each map $T\mapsto Tb_k$ is continuous on $\mathcal{F}(E\oplus C_1)$, $Y$ is a closed subspace of $\mathcal{F}(E\oplus C_1)^{\oplus n}$. It is clear, furthermore, that $Y$ is a left $\mathcal{F}(E\oplus C_1)$-submodule of $\mathcal{F}(E\oplus C_1)^{\oplus n}$.
\paragraph*{}Since $E\oplus C_1$ has the approximation property, for each $\epsilon >0$ there is a finite rank operator, $U$, in $\mathcal{F}(E\oplus C_1)$ such that $\sum_{k=1}^n ||UT_k-T_k|| < \epsilon$, where $(UT_1,UT_2,...\, , UT_n)$ belongs to $Y$ because $Y$ is a left $\mathcal{F}(E\oplus C_1)$-module. Therefore for each k there is a sequence $\{ T_k^{(p)} \}_{p=1}^\infty$ of finite rank operators on $E\oplus C_1$ such that:
\\
(i) for each $k$, $\sum_{p=1}^\infty ||T_k^{(p)}|| < \infty$ and $\sum_{p=1}^\infty T_k^{(p)} = T_k $; and
\\
(ii) for each $p$, $(T_1^{(p)}, T_2^{(p)}, ... \, ,T_n^{(p)})$ belongs to $Y$. 

\paragraph*{}For each $p$ choose a finite dimensional subspace, $E_p$, of $E\oplus C_1$ such that the range of $T_k^{(p)}$ is contained in $E_p$ for each $k$. Then choose a finite dimensional space $G_{i_p}$ such that $d(E_p, G_{i_p})< 2$ and an isomorphism $U_p : E_p \to G_{i_p}$ such that $||U_p|| =1$ and $||U_p^{-1}||<2$. Finally, choose a sequence of positive numbers $\{\lambda_p \}_{p=1}^\infty$ which converges to zero and is such that $\sum_{p=1}^\infty ||\lambda^{-1} T_k^{(p)} || <\infty$ for each $k$. The spaces $G_{i_p}$ are subspaces of $C_1$, and hence of $E\oplus C_1$, and there is a projection, $Q_p$, from $E\oplus C_1$ onto $G_{i_p}$ with $||Q_p||=1$. Define operators $R_k, \, k=1,2,...\, , n$ and $V$ in $\mathcal{F}(E\oplus C_1 )$ by
$$ R_k = \sum_{p=1}^\infty \lambda_p^{-1} U_p T_k^{(p)} \quad \text{and} \quad V = \sum_{p=1}^\infty \lambda_p U_p^{-1} Q_p. $$ 

\paragraph*{}Let, for each $P$, $R_k^{(P)} = \sum_{p=1}^P \lambda_p^{-1}U_pT_k^{(p)}$. Then, since $Y$ is a left $\mathcal{F}(E\oplus C_1)$-module,  $(R_1^{(P)}, R_2^{(P)}, ... \, , R_n^{(P)})$ belongs to $Y$ for each $P$. Since $Y$ is closed, it follows that $(R_1, R_2, ...\, , R_n)$ belongs to $Y$, that is, that $\sum_{k=1}^n Y_k b_k =0$. Also, since $Q_pU_q =0$ if $p\neq q$, we have $VR_k = \sum_{p=1}^\infty U_p ^{-1}Q_pU_pT_k^{(p)} = \sum_{p=1}^\infty T_k^{(p)} = T_k$. Therefore $\mathcal{F}(E\oplus C_1)$ satisfies the characterization of right universal flatness quoted above.  
\paragraph*{}It is clear that (3) implies (2). To show that (2) implies (1), suppose that $E$ does not have the approximation property. By a theorem of Grothendieck, see [LT, theorem 1.e.4], there is a Banach space, $F$, a compact operator $T: F \to E\oplus C_1$ such that $T$ cannot be approximated in norm by finite rank operators. It may be seen from the proof of the theorem that $F$ may be chosen to separable and thus be a quotient space of $\ell_1$. Now $E\oplus C_1 \simeq G\oplus \ell_1$, for some space $G$, and so the quotient map may be extended to yield a surjection $Q: E\oplus C_1 \to F$ by defining the extension to zero on $G$. Put $R= TQ$. Then $R$ belongs to $\mathcal{F}(E\oplus C_1)$ because it is compact and effectively an operator from $\ell_1$. Let S be any compact operator on $\ell_1$ whose range is a dense subspace of the kernel of the quotient map $\ell_1 \to F$ and extend $S$ to be an operator on $E\oplus C_1$ by defining its extension to be zero on $G$. Then $S$ belongs to $\mathcal{F}(E\oplus C_1)$ and $RS =0$.
\paragraph*{}If $\mathcal{F}(E\oplus C_1)$ were right flat there would be $A_1, A_2, ... \, , A_n$ and $B_1,B_2, ... \, , B_n$ in $\mathcal{F}(E\oplus C_1)$ such that $R=\sum_{k=1}^n A_kB_k$ and $B_kS=0$ for each $k$. The condition $B_kS =0$ implies that  the kernel of $B_k$ would contain the kernel of the quotient map $\ell_1 \to F$ and so $B_k$ would induce a map $B'_k : F \to E\oplus C_1$. We would then have that $T = \sum_{k=1}^n A_kB'_k$ and the requirement that each $A_k$ belong to $\mathcal{F}(E\oplus C_1)$ would imply that $T$ belongs to $\mathcal{F}(F, E\oplus C_1)$. This would contradict the choice of $T$ and so $\mathcal{F}(E\oplus C_1)$ is not right flat.
\paragraph*{}The reason that $C_1$ appears in the theorem, rather than $C_p$ for general $p$, is to guarantee in the second part of the proof that $\ell_1$ is a direct summand in $E\oplus C_1$. Thus $E$ having the approximation property is also equivalent to the right (universal) flatness of $\mathcal{F}(F, E\oplus C_p\oplus \ell_p)$. We also have the

\begin{corollary}
Let $E$ be a Banach space with the approximation property. Then $\mathcal{F}(E\oplus C_p)$ is right universally flat.
\end{corollary}

\section*{Left Universal Flatness of $\mathcal{F}(E)$}
\paragraph*{}An argument dual to the first part of the proof of theorem 1 proves the next

\begin{corollary}
Let $E$ be a Banach space such that $E^*$ has the approximation property. Then $\mathcal{F}(E\oplus C_p)$ is left universally flat.
\end{corollary}   

\paragraph*{}The approximation property for the dual space may also be characterized in terms of a homological property. Spaces whose duals do not have the approximation property are characterized in theorem 1.e.5 from [LT]. This theorem says that $E^*$ does not have the approximation property if and only if there is a space $F$ and a compact operator $T:E\to F$ which cannot be approximated by finite rank operators. The theorem may be used in a dual way to that in which theorem 1.e.4 was used above. However, to carry this through it is necessary to embed $F$ in some space which has the approximation property. None of the spaces $C_p$ have the property that any separable space can be embedded in them. The space $C[0,1]$ of continuous functions on the unit interval does have this property, see [Woj, theorem II.B.4], and so, as a dual to theorem 1, we have

\begin{thm}
For each Banach space $E$ the following are equivalent
\\
\\
$(1)$ $E^*$ has the approximation property;
\\
\\
$(2)$ $\mathcal{F}(E\oplus C_p \oplus C[0,1])$ is left flat; \quad and
\\
\\
$(3)$ $\mathcal{F}(E\oplus C_p \oplus C[0,1])$ is left universally flat.
 
\end{thm}
 
\paragraph*{}If $E^*$ has the approximation property, then $E$ has the approximation property, see [LT, theorem 1.e.7]. Hence it follows from theorems 1 and 2 that, if $\mathcal{F}(E\oplus C_1 \oplus C[0,1])$ is left flat, then it is right universally flat. There are Banach spaces, $E$, such that $E$ has the approximation property but $E^*$ does not, see [LT, theorem 1.e.7]. It follows that the right flatness of $\mathcal{F}(E\oplus C_1 \oplus C[0,1])$ does not imply its left flatness.

\paragraph*{}If $\mathcal{K}(E)$ has a bounded right approximation identity, then it is left universally flat. It is shown in [GW, corollary 2.7] that, if $\mathcal{K}(E)$ has a bounded right approximate identity, then it has a bounded left approximate identity and so it will also be right universally flat. There are spaces for which $\mathcal{K}(E)$ has a bounded left approximate identity but not a bounded right identity.

\section*{Remarks and Open Questions}
Flatness of algebras was proved in [Wo3] by means of a factorization property ($\Phi$) and H-unitality was proved in [Wo2] by means of another factorization property (F). (Note that ($\Phi$) is a purely algebraic condition whereas (F) applies to topological algebras.) If $E$ has the bounded approximation property, then $\mathcal{F}(E)$ has a bounded left approximate identity and it follows that both ($\Phi$) and (F) hold in $\mathcal{F}(E)$. The above arguments do not establish that either ($\Phi$) or (F) holds in $\mathcal{F}(E\oplus C_1)$ when $E$ has the approximation property. 

\newtheorem{ques}{Question}

\begin{ques}
If $E$ has the approximation property must $\mathcal{F}(E\oplus C_1)$ satisfy either $(\Phi)$ or $(F)$?
\end{ques}
In the proofs of the equivalences in theorem 1, the direct summand $C_1$ was required in both directions. These results say nothing about $\mathcal{F}(E)$ when $E$ does not have a direct summand isomorphic to $C_1$.

\begin{ques}
Is $\mathcal{F}(E)$ right flat or right universally flat whenever $E$ has the approximation property? 
\end{ques}

\begin{ques}
Does $E$ have the approximation property whenever $\mathcal{F}(E)$ is right flat or right universally flat?
\end{ques}
Theorem 1 adds to the class of Banach spaces for it is known that $\mathcal{K}(E)$ is H-unital. It also shows that there are Banach spaces that are not right flat. However, there are no spaces known to the author for which $\mathcal{K}(E)$ is not H-unital.

\begin{ques}
Are $\mathcal{K}(E)$ and $\mathcal{F}(E)$ always H-unital?
\end{ques}
A special case of this question is whether the sequence (1) is exact at $A$ for these algebras.

\begin{ques}
Is there a Banach space such that the linear span of all products in $\mathcal{F}(E)$ or $\mathcal{K}(E)$ is not equal to $\mathcal{F}(E)$ or $\mathcal{K}(E)$ respectively?
\end{ques}
Theorem 1 in [Jo] shows that $\mathcal{F}(E\oplus C_p)$ is equal to the span of its products for any Banach space $E$. This makes the question as to whether $\mathcal{F}(E\oplus C_1)$ is H-unital particularly interesting.
\paragraph*{}The idea behind theorem 1 in [Jo] may be used to prove much more than factorization of single elements. We have the following

\newtheorem*{prop}{Proposition}

\begin{prop}
Let $E$ be a Banach space and $\{T_n \}_{n=1}^\infty$ be a sequence of operators in $\mathcal{F}(E\oplus C_p)$ such that $||T_n|| \to 0$ as $n \to \infty$. Then there are:
\\
\\
$(1)$ an element $U$ and a sequence $\{S_n \}_{n=1}^\infty$ in $\mathcal{F}(E\oplus C_p)$ such that $T_n =US_n$ for each $n$ and $||S_n||\to 0$ as $n\to \infty$; \, and
\\
\\
$(2)$ an element $U'$ and a sequence $\{S'_n \}_{n=1}^\infty$ in $\mathcal{F}(E\oplus C_p)$ such that $T_n =S'_nU'$ for each $n$ and $||S'_n||\to 0$ as $n\to \infty$.
\end{prop}

\paragraph*{}The proof of (1) is essentially the same argument as used in [Jo] while the proof of (2) is just the dual of that argument. In general, $\mathcal{F}(E\oplus C_p)$ does not have a right or left bounded approximate identity and so the proposition provides more examples of Banach algebras which satisfy a strong factorization property without having a bounded approximate identity, see [DW] and [Wi2].
\paragraph*{}The next result is an immediate consequence of this proposition.

\begin{corollary}
\label{cor:four}
Let $E$ be a Banach space and $T$ be either a right or left Banach $\mathcal{F}(E\oplus C_p)$-module homomorphism from $\mathcal{F}(E\oplus C_p)$. Then $T$ is continuous.
\end{corollary}
\paragraph*{}A right multiplier on an algebra is just a left module homomorphism from the algebra to itself. The corollary thus provides another proof of the lemma for spaces which have $C_p$ as a direct summand. It was to such spaces that we applied the lemma and so we might just as well have used this proof.
\pagebreak

\end{document}